\documentclass{amsart}

\usepackage{amscd}

\newtheorem{theorem}{Theorem}[section]

\newtheorem{proposition}[theorem]{Proposition}

\theoremstyle{definition}     % italic or bold etc.

\theoremstyle{remark}

\numberwithin{equation}{section}

\newcommand{\Pmatrix}{\left( \begin{matrix} }
\newcommand{\Endpmatrix}{\end{matrix} \right) }

\begin{document}

\title[Kummer structures on a K3 surface]{Kummer structures on a K3 surface
- An old question of T. Shioda}

\author[S. Hosono]{Shinobu Hosono}

\address{
Department of Mathematical Sciences,
University of Tokyo, Komaba Meguro-ku,
Tokyo 153-8914, Japan
}
\email{hosono@ms.u-tokyo.ac.jp}

\author[B.H. Lian]{Bong H. Lian}

\address{
Department of Mathematics, Brandeis University,
Waltham, MA 02154, U.S.A.}
\email{lian@brandeis.edu}

\author[K. Oguiso]{Keiji Oguiso}

\address{
Department of Mathematical Sciences,
University of Tokyo, Komaba Meguro-ku,
Tokyo 153-8914, Japan
}
\email{oguiso@ms.u-tokyo.ac.jp}

\author[S.-T. Yau]{Shing-Tung Yau}

\address{
Department of Mathematics, Harvard University,
Cambridge, MA 02138, U.S.A.}
\email{yau@math.harvard.edu}

\dedicatory{Dedicated to Professor Tetsuji Shioda on the occasion of
his sixtieth birthday}


\begin{abstract}
We apply our earlier results on Fourier-Mukai partners
to answer definitively a question about
Kummer surface structures, posed by T. Shioda 25 years ago.
\end{abstract}

\maketitle

\noindent
\centerline{0. {\bf Introduction} }
\par
\vskip 4pt

In this note, unless stated otherwise, we shall work in the
category of smooth complex projective varieties. For a smooth
projective variety $X$, we denote by $D(X)$ the bounded derived
category of coherent sheaves on $X$. The object $Ob D(X)$ contains
the set $\{\mathcal O_{x} \vert x \in X\}$ of structure sheaves of
closed points of $X$. We regard $D(X)$ as a triangulated category
in a natural manner. By an equivalence $D(Y) \simeq D(X)$ we
always mean an equivalence of triangulated categories.
\par
\vskip 4pt
In [BO], Bondal and Orlov showed that a lot of information
can be extracted from the objects
$\{\mathcal O_{x} \vert x \in X\}$ when the canonical divisor
$K_{X}$ is ample or anti-ample. In particular, in this case they further showed that
$Y \simeq X$ if $D(Y) \simeq D(X)$.
Quite recently, Kawamata [Ka] obtained a generalization of this result. He showed
that $Y$ is birationally
equivalent to $X$ if $D(X) \simeq D(Y)$ and if $K_{X}$ or $-K_{X}$ is big.
\par
\vskip 4pt
On the other hand, Mukai [Mu1] showed, about 20 years ago, that the
Poincar\'e
bundle $\mathcal P$ on $A \times \hat A$, where $A$ is an abelian variety
and $\hat A := \text{Pic}^{0} A$, induces an equivalence
$\Phi^{\mathcal P}_{\hat A \rightarrow A}: D(\hat A) \rightarrow D(A)$.
The functor $\Phi^{\mathcal P}_{\hat A \rightarrow A}$ is the so-called
Fourier-Mukai transform defined by
\begin{equation}
\Phi^{\mathcal P}_{\hat A \rightarrow A}(\mathcal X)=
\mathbf R \pi_{A*}( \mathcal P \otimes^{\hskip-7pt \,^\mathbf L}
\mathbf L \pi_{\hat A}^*
\mathcal X )
\label{eqn:FM}
\end{equation}
where $\pi_{\hat A}: \hat A \times A \rightarrow \hat A$ and
$\pi_A : \hat A \times A \rightarrow A$ are the natural projections. In this
equivalence, the structure sheaf $\mathcal O_{\hat a}$ of the point
$\hat a \in \hat A$
is mapped to the invertible sheaf $\mathcal P_{\hat a}$ on $A$
corresponding to $\hat a$. Therefore the derived category of an abelian
variety does not characterize
the structure sheaves of points anymore. Indeed, $\hat A$ is not
(even birationally) isomorphic to $A$ in general even when $D(\hat A) \simeq
D(A)$.
\par
\vskip 4pt
This example suggests that when $K_X$ is trivial,
$D(X)$ no longer has enough information to reconstruct
(the birationally equivalence class of) the variety $X$,
and that a new and interesting relationship
arises between varieties $X,Y$ with $D(X)\simeq D(Y)$.
\par
\vskip 4pt
The aim of this note is to examplify
this idea by studying one old problem in concrete geometry, namely
the following
problem posed by T. Shioda 25 years ago, from the view point of Fourier-Mukai
partners:
\par
\vskip 4pt
\noindent
{\bf Problem.} ([Sh1, Question 5]) Does the Kummer variety
$\text{Km}\, A$ uniquely determine
the abelian variety $A$ up to isomorphism, i.e.
$\text{Km}\, B \simeq \text{Km}\, A \,
\Rightarrow \, B \simeq A$?
\par
\vskip 4pt
As shown by Shioda, the answer is affirmative in dimension $\geq 3$,
and it is also affirmative if $\text{dim}\, \text{Km}\, A = 2$ and
$\rho(\text{Km}\, A) = 20$, i.e. $\rho(A) = 4$,
the maximal possible case [Sh1],[SM, Theorem 5.1]. (See also Section 2.)
He then expected an affirmative answer
to the Problem in dimension two.
\par
\vskip 4pt
However, as it is expected from Shioda's Torelli Theorem [Sh2] for
abelian surfaces, which he found after he posed the problem above, 
that one has
$\text{Km}\, \hat A \simeq \text{Km}\, A$.
This was first noticed by [GH, Theorem 1.5 and Remark in Sect.I].
See also [HS, Sect.III.3], and Section 1 for another explanation.
Therefore, it is natural and interesting to
consider the following modified:
\par
\vskip 4pt
\noindent
{\bf Problem.} Let $A$ be an abelian surface.
Then, $\text{Km}\, B \simeq \text{Km}\, A \,
\Rightarrow \, B \simeq A \, \text{or} \, B \simeq \hat A$?
\par
\vskip 4pt

For our statement, we need a few preparations. First, we set
$$
FM(X) := \{Y \vert D(Y) \simeq D(X)\}/\text{isom}\, .$$
An element of $FM(X)$ is called a Fourier-Mukai (FM) partners of $X$.
Recall that if $X$
is an abelian (resp. K3) surface then so are its FM partners $Y$, and that
$$
FM(X) = \{Y \vert (T(Y), \mathbf C \omega_{Y})\, \simeq\,
(T(X), \mathbf C \omega_{X})\, \}/
\text{isom}\, .
$$
Here $T(X)$ denotes the transcendental
lattice of $X$, i.e. $T(X) := NS(X)^{\perp}$ in $H^{2}(X, \mathbf Z)$,
$\omega_X$ denotes a nowhere vanishing holomorphic two form on $X$,
and the isomorphism
$(T(Y), \mathbf C \omega_{Y})\, \simeq\, (T(X), \mathbf C \omega_{X})$ stands
for a Hodge isometry.
These properties are due to Mukai and Orlov ([Mu2], [Or]) and are also well
treated in [BM, Theorem 5.1].
\par
\vskip 4pt
Next, for a K3 surface $X$, we define
$$
\mathcal K(X) := \{\, B \; \vert \;
B\, \text{is an abelian surface s.t.}\, \text{Km}\,
B \simeq\, X \, \}/\text{isom}\, .
$$
We call $\mathcal K(X)$ the set of Kummer structures on $X$. This set
$\mathcal K(X)$ measures how many different Kummer surface structures $X$ has.
Note that $\mathcal K(X) = \emptyset$ unless $X$ is a Kummer surface.
\par
\vskip 4pt

Our main result is now stated as follows:

\begin{theorem} {\bf (Main Theorem)} Let $A$ be an abelian surface and
assume that $X = \text{Km}\, A$. Then:
%%%%%%%
\begin{enumerate}
\item $\mathcal K(X) = FM(A)$. In particular, $\vert \mathcal K(X)
\vert < \infty$ and $\hat A \in \mathcal K(X)$.
%%%%%%%
\item If $\rho(A) = 3$ and $\text{det}\, NS(A)$ is square free or
if $\rho(A) = 4$,
then
$\mathcal K(X) = \{A, \hat A\}$.
%%%%%%%%
\item For an aribitrarily given natural number $N$, there is an abelian
surface $A$ such that $\mathcal K(X)$, where $X = \text{Km}\, A$, contains
$N$ abelian surfaces
$A_{i}$ ($1 \le i \le N$) with $NS(A_{i}) \not\simeq NS(A_{j})$, hence
$A_{i} \not\simeq A_{j}$, for all $1 \le i\not= j \le N$.
\end{enumerate}
\end{theorem}

Assertion 1 makes an interesting connection between two important
notions: Kummer surface structures and Fourier-Mukai partners.
Assertion 2 is probably not surprising in light of Shioda's
earlier result. However,
Assertion 3 shows that
K3 surfaces can have arbitrarily high number of non-isomorphic
Kummer surface structures! This is quite unexpected
since it runs counter to Shioda's original observation.
It is also surprising in the light of the theory
of lattices. (See also Theorem 3.1 and Remark in \S 3.) Assertions 1, 2,
and 3 will be proved respectively in Sections 1, 2, and 3.
\par
\vskip 4pt
In Appendix, we remark that $\hat A \not\in \mathcal K(\text{Km}\, A)$ while
$\hat A \in FM(A)$, whence $\mathcal K(\text{Km}\, A) \not= FM(A)$
in general, if a complex $2$-torus $A$ is not projective. This is pointed out
to us by Yoshinori Namikawa. See also [Na] for relevant work on generalized
Kummer varieties and
differences between dimension $2$ and dimension $\ge 4$.
\par
\vskip 4pt
Throughout this note, an abelian surface means a projective complex $2$-torus,
while a complex $2$-torus is not assumed to be projective.
\par
\vskip 4pt
\noindent
{\it Acknowledgement.} We thank Professors
B. Gross, Y. Namikawa and T. Shioda for their interests in this work
and several valuable comments. The main part of this work has been
done during the first and third named authors' stay
at Harvard University. They would like to express their thanks to the
Harvard University and the Education Ministry of Japan
for financial support. The second named author is supported by
NSF DMS-0072158.

\section{\bf FM partners and Kummer structures}

In this section, we shall prove assertion 1 of Main Theorem.
\par
\vskip 4pt
\noindent
{\it Proof of assertion 1.} Let $A$ and $B$ be abelian surfaces.
By [Ni1] (see also [Mo, Lemma 3.1]), the canonical rational map
$\pi_{A} : A  \; \text{-}\,\text{-}\,\text{-} \hskip-6pt >  \text{Km}\, A$
gives the Hodge isometry
$$
\pi_{A, *} : (T(A)(2), \mathbf C \omega_{A}) \simeq
(T(\text{Km}\, A), \mathbf C \omega_{\text{Km}\, A})\, ,
$$
and likewise for $B$ and $\text{Km}\, B$. Here, for a lattice
$L := (L, (*, **))$ and
for a non-zero integer $m$,
we define the lattice $L(m)$ by $L(m) := (L, m(*, **))$.
\par
\vskip 4pt
Therefore, one has
$$
(T(B), \mathbf C \omega_{B}) \simeq (T(A), \mathbf C \omega_{A})\,
\Longleftrightarrow\, (T(\text{Km}\, B), \mathbf C \omega_{\text{Km}\, B})
\simeq
(T(\text{Km}\, A), \mathbf C \omega_{\text{Km}\, A})\, .
$$
Hence, by the characterization of the FM partners of K3 and abelian surfaces,
one obtains
$$B \in FM(A)\, \Longleftrightarrow\, \text{Km}\, B \in FM(\text{Km}\, A)\, .$$
On the other hand, we know that $FM(\text{Km} A) = \{\text{Km}\, A\}$
by [Mu2] or by
the Counting Formula [HLOY, \S 3]. Then we have
$$\text{Km}\, B \in FM(\text{Km}\, A)\, \Longleftrightarrow\, \text{Km}\, B
\simeq \text{Km}\, A\, ,$$
and therefore
$$B \in FM(A)\, \Longleftrightarrow\, B \in \mathcal K(X)\, .$$
\par
\vskip 4pt
Since $\vert FM(A) \vert < \infty$ by [BM, Proposition 5.3] and
$\mathcal K(X) = FM(A)$,
we have now $\vert \mathcal K(A) \vert < \infty$. See also Section 2 for a
more explicit estimate
of $\vert FM(A) \vert$.  As it is remarked in
Introduction, $\hat A \in FM(A)$. Thus, by $\mathcal K(X) = FM(A)$, we have
$\hat A \in \mathcal K(X)$ as well. \qed

\section{\bf Kummer structures in larger Picard numbers}

In this section, we shall prove assertion 2 of Main Theorem.
For this purpose, we shall first give an effective estimate
of $\vert FM(A) \vert$. This is based on the argument of
[HLOY, \S 3] and the following Theorem due to Shioda [Sh2]
(See also [Mo]):

\begin{theorem}\label{thm:Shioda}
\par
\begin{enumerate}
\item Let $A$ and $B$ be complex $2$-tori. Then
$$(H^{2}(B, \mathbf Z), \mathbf C \omega_{B}) \simeq
(H^{2}(A, \mathbf Z), \mathbf C \omega_{A})\, \Longleftrightarrow\,
B \simeq A\, {\rm or}\, B \simeq \hat A\, .$$
\item For any given weight two Hodge structure
$(U^{\oplus 3}, \mathbf C \omega)$, there is a marked complex $2$-torus
$(A, \tau_{A})$ such that $\tau_{A}$ gives a Hodge isometry
$$
\tau_{A} : (H^{2}(A, \mathbf Z), \mathbf C \omega_{A})
\simeq (U^{\oplus 3}, \mathbf C \omega)\, .
$$

\end{enumerate}
\end{theorem}
For the second statement, we recall that $U$ is the even unimodular
hyperbolic lattice of rank $2$ and the second
cohomology lattice $H^{2}(A, \mathbf Z)$ of a complex $2$-torus $A$
is always isometric to the lattice $U^{\oplus 3}$.
\par
\vskip 4pt

Let $A$ be an abelian surface. Denote by
$G = O_{Hodge}(T(A),\mathbf C \omega_{A})$ the group of Hodge
isometries of $(T(A), \mathbf C \omega_{A})$. We call two primitive
embeddings $\iota:
T(A) \rightarrow U^{\oplus 3}$, $\iota': T(A) \rightarrow U^{\oplus 3}$
$G$-{\it equivalent} if there exist $\Phi \in O(U^{\oplus 3})$ and $g \in G$
such that $\iota' \circ g = \Phi \circ \iota$, i.e. the following diagram
commutes;
$$
\begin{CD}
\iota: T(A) @>>> U^{\oplus 3} \cr
  @VV{\exists g}V    @VV{\exists \Phi}V \cr
\iota': T(A) @>>> U^{\oplus 3}  \cr
\end{CD}
$$
We denote by $\mathcal P^{G\text{-eq}}(T(A),U^{\oplus 3})$ the set of
$G$-equivalence classes
of primitive embeddings of $T(A)$ into $U^{\oplus 3}$:
$$
\mathcal P^{G\text{-eq}}(T(A),U^{\oplus 3}):=
\{ \text{primitive embedding } \iota: T(A) \rightarrow U^{\oplus 3} \}/
G\text{-equiv.} \;\;.
$$

\vskip0.3cm

\noindent
{\bf Claim.}
{\it There is a surjection,}
$$
\xi: FM(A) \longrightarrow \mathcal P^{G\text{-eq}}(T(A),U^{\oplus 3})\,,
\; B \mapsto [\iota_B] \;,
$$
{\it such that} $\xi^{-1}([\iota_B])=\{ B , \hat B \}$.

\begin{proof} Let $B \in FM(A)$. Choose a marking $\tau_B: H^2(B,\mathbf Z)
\rightarrow U^{\oplus 3}$ and a Hodge isometry $\gamma:
(T(A), \mathbf C \omega_A) \simeq (T(B), \mathbf C \omega_B)$.
Define $\xi: [B] \mapsto [\iota_B]$, where $\iota_B=\tau_B \circ \gamma$.
This is a well-defined surjection with $\xi^{-1}(\xi(A)) = \{A, \hat A\}$.
The proof of this fact is
exactly the same as the argument
for Theorem 3.3 of [HLOY], except that Theorem 2.1 here plays
the corresponding role of the Torelli Theorem for K3 surfaces there.
Note that $\xi$ may fail to be injective because $B$ and $\hat B$ need not
be isomorphic in general.
\end{proof}

Recall Theorem 1.4 of [HLOY]:
\begin{equation}
\vert \mathcal P^{G\text{-eq}}(T(A), U^{\oplus 3}) \vert
=  \sum_{j=1}^m |O(S_j)\setminus O(A_{S_j}) / G| \;\; .
\label{eqn:CFM}
\end{equation}
Here $\mathcal G(NS(A)) :=\{ S_1,S_2,\cdots , S_m \} \; (S_1 \simeq NS(A))$
is the genus of the N\'eron-Severi lattice
$NS(A)$ and $O(A_{S_{i}})$ is the orthogonal group of the
discriminant group $A_{S_{i}} = S_{i}^{*}/S_{i}$ with respect to the natural
$\mathbf Q/2 \mathbf Z$-valued quadratic form $q_{S_{i}}$ and the set
$O(S_j)\setminus O(A_{S_j}) / G$
is the orbit space of the natural action of $O(S_j) \times G$ on $A_{S_{j}}$.
Note that
$NS(A)$ is even and hyperbolic. We also recall that the genus of $NS(A)$ is
nothing but the set of
isomorphism classes
of even lattices $S$ such that $(NS(A), q_{NS(A)}) \simeq (S, q_{S})$,
by [Ni2, Corollary 1.9.4].
\par
\vskip 4pt

We can now prove assertion 2 of Main Theorem.
\par
\vskip 4pt
If $\text{det}\,NS(A)$ is
square free, it follows
that $A_{NS(A)}$ is cyclic, i.e. $l(A_{NS(A)}) = 1$. Here $l(M)$ denote
the minimal number of generators of a finite abelian group $M$. Hence,
if, in addition, $\rho(A) = 3$, then one has
$$\text{rank}\;NS(A)=3 \geq 2+l(A_{NS(A)})\, .$$
Similarly, if $\rho(A) = 4$, then one has
$$\text{rank}\;NS(A) = 4 \geq 2+l(A_{NS(A)})\, .$$
Indeed, if $\rho(A) = 4$, then $\text{rank}\, T(A) = 2$.
Since $(A_{NS(A)}, q_{NS(A)}) \simeq
(A_{T(A)}, -q_{T(A)})$, we have
$l(A_{NS(A)}) = l(A_{T(A)}) \le \text{rank}\, T(A) = 2$.
\par
\vskip 4pt
Hence, by [Ni2, Theorem 1.14.2], one has that $\mathcal G(NS(A)) =
\{NS(A)\}$ and that the natural map $O(NS(A)) \rightarrow O(A_{NS(A)})$ is
surjective in each case. Now it follows from the equation (2.1) that
$\vert \mathcal P^{G\text{-eq}}(T(A), U^{\oplus 3}) \vert =1$.
Hence, by the preceding claim, we have $FM(A)=\{ A, \hat A\}$.
Combining this with assertion 1 of Main Theorem, we obtain
the desired equality $\mathcal K(X) = \{A, \hat A\}$.
\qed
\par
\vskip 4pt
\noindent
{\bf Remark.} In [SM], it is shown that an abelian surface $A$ with
$\rho(A) = 4$ is necessarily a direct product of two elliptic curves, namely,
$A \simeq E \times F$. Since $\hat E \simeq E$ and $\hat F \simeq F$,
it follows that $\hat A \simeq A$ if $A \simeq E \times F$. Combining this
with assertion 2, one has $\mathcal K(X) = \{A \}$ if $\rho(A) = 4$.
This gives an alternative explanation for the result [SM, Theorem 5.1]
quoted in Introduction.
\qed
\section{\bf K3 surfaces with many Kummer structures}

In this section, we shall show assertion 3 of Main Theorem.
\par
\vskip 4pt
We begin by the following result first found in [Og2]:
\begin{theorem}\label{thm:ogiw} For any given natural number $N$,
there is an even hyperbolic lattice $S$ such that $\text{rank}\,
S = 2$, $\text{det}\, S = -pq$
and $\vert \mathcal G(S) \vert \geq N$.
Here $\mathcal G(S)$ is the genus of the lattice $S$, and $p$ and $q$
are some distinct
prime numbers (including $1$).
\end{theorem}
\par
\vskip 4pt
\noindent
{\bf Remark.} Such phenomena seem quite rare and are perhaps impossible if
$\text{rank}\, S \geq 3$ (cf. [Cs, Chapters 9 and 10] and
[CL, \S 8 - \S10]). In this sense,
as it is pointed out to us by Professor B. Gross, this Theorem is a sort
of miracle which happens probably only in rank $2$. \qed

\begin{proof} Here we shall sketch the argument. Our argument requires
some basic
properties about class groups. These are found  in the book [Za].
(See also [HLOY, Section 4]
for summary.) Let $P$ be the
set of prime numbers including $1$.  By [Iw, Main Theorem],
the set
$$\mathcal N := \{(n, p, q) \in \mathbf N \times P^{2} \vert
4n^2 + 1 = pq\}$$
is an infinite set. Note that $pq \equiv 1\, \text{mod}\, 4$ and $p \not= q$
for $(n, p, q) \in \mathcal N$. Therefore $pq$ is a so-called fundamental
discriminant. Let $H(pq)$ be the narrow class group of (the ring of integers
$O(pq)$ of) the real quadratic field $\mathbf Q(\sqrt{pq})$.
We denote by $h(pq)$ the class
number of $\mathbf Q(\sqrt{pq})$, i.e. $\vert H(pq) \vert$. Recall that
there is a natural bijective correspondence between $H(pq)$ and the set
$\mathcal L(pq)$ of the proper, i.e. $\text{SL}(2, \mathbf Z)$, isomorhism
classes of even hyperbolic rank $2$ lattices $S$ with $det\, S = -pq$. Recall
also that $\mathcal L(pq)$ is decomposed into at most two
sets of the same cardinality:
$$\mathcal L(pq) = \tilde{\mathcal G}(S_1) \cup \tilde{\mathcal G}(S_2)\, .$$
Here $\tilde{\mathcal G}(S_i)$ is the set of proper, i.e.
$\text{SL}(2, \mathbf Z)$, isomorphism classes of lattices in the same genus as
$S_{i}$. Combining these two facts with
$\text{GL}(2, \mathbf Z)/\text{SL}(2, \mathbf Z) \simeq \mathbf Z/2$, we have
$$\vert \mathcal G(S_{1})
\vert \ge \frac{\vert \tilde{\mathcal G}(S_{1}) \vert}{2} \ge
\frac{\vert \mathcal L(pq) \vert}{4} = \frac{h(pq)}{4}\, .$$
Therefore, it is enough to show that $h(pq) \rightarrow \infty$ when
$pq \rightarrow \infty$. We note that there is a sequence
$\{(n_{k}, p_{k}, q_{k})\}_{k=1}^{\infty} \in \mathcal N$
such that $p_{k}q_{k} \rightarrow \infty$ when $k \rightarrow \infty$. This
is because $\vert \mathcal N \vert = \infty$.
\par
\vskip 4pt
Let us show that $h(pq) \rightarrow \infty$ when
$pq \rightarrow \infty$. Since  $(2n + \sqrt{pq})(-2n + \sqrt{pq}) = 1$
by $4n^2 + 1 = pq$,
the element $2n + \sqrt{pq} > 1$ is a unit of $O(pq)$. Therefore, the
fundamental unit $\epsilon(pq)$ of $O(pq)$ satisfies $1 < \epsilon(pq) < pq$.
Now, combining this with the Siegel-Brauer formula
$$\lim_{pq \rightarrow \infty}
\frac{\log (h(pq)\log \epsilon (pq))}{\log pq} = \frac{1}{2}\, ,$$
one has $h(pq) \rightarrow \infty$ when $pq \rightarrow \infty$.
\end{proof}

Let us return to the proof of assertion 3 of Main Theorem.
\par
\vskip 4pt
Let $S$ be the lattice in Theorem (4.1) and take $N$ non-isomorphic
elements
$$
S_{1} := S,  S_{2}, \cdots , S_{N} \in \mathcal G(S) \;.
$$
Then $S_{i}$ are even hyperbolic lattices of $\text{rank}\,
S_{i} = 2$ and satisfy $(A_{S_{i}}, q_{S_{i}}) \simeq (A_{S}, q_{S})$.
\par
\vskip 4pt
Let us choose a primitive embedding
$$
\varphi_{i} : S_{i} \rightarrow U^{\oplus 3}
$$
for each $i$. For instance, if
$S_{i} = \mathbf Z\langle v_{i1}, v_{i2} \rangle$ with
$$((v_{i1}, v_{i2})) =
\Pmatrix  2a_{i} & b_{i}\\ b_{i} & 2c_{i}\\ \Endpmatrix ,$$
then we can define $\varphi_{i}$ by
$\varphi_{i}(v_{i1}) = e_{1} + a_{i}f_{1}$ and
$\varphi_{i}(v_{i2}) = b_{i}f_{1} + e_{2} + c_{i}f_{2}$, where $e_{j}$
and $f_{j}$ are the standard basis of
the $j$-th factor $U$ of $U^{\oplus 3}$.
\par
\vskip 4pt
Set $T_{i} := \varphi_{i}(S_{i})^{\perp}$ in $U^{\oplus 3}$
and $T := T_{1}$. Then
$$
(A_{T_{i}}, q_{T_{i}}) \simeq (A_{S_{i}}, -q_{S_{i}}) \simeq (A_{S}, -q_{S})
$$
and the signature of $T_{i}$ is $(2, 2)$ for each $i$. Therefore
$T_{i} \in \mathcal G(T)$.
\par
\vskip 4pt
Since $pq$ is square free, we have $A_{T} \simeq A_{S} \simeq \mathbf Z/pq$.
Thus,
$$
\text{rank}\, T = 4 \geq 2 + 1 = 2 + l(A_{T}) \; .
$$
In addition, $T$ is indefinite. Hence $\mathcal G(T) = \{T\}$
by [Ni2, Theorem 1.14.2].
Therefore $T_{i} \simeq T$ for each $i$. Let us choose an isometory
$$
\sigma_{i} : T \; \simeq \; T_{i}
$$
for each $i$. We set $\sigma_{1} := id$.
\par
\vskip 4pt
Choose a maximal weight two
Hodge structure $(T, \mathbf C \omega)$ on $T$. Here
the term {\it maximal} means that $T$ satisfies the following
property: if $T' \subset T$ is a primitive sublattice of $T$ such that
$\mathbf C \omega \in T' \otimes \mathbf C$ then $T' = T$. This is
equivalent to that $\langle \text{Re} \, \omega, \text{Im} \, \omega
\rangle^{\perp}
\cap T = \{0\}$ in $T \otimes \mathbf R$. As it is easily shown,
there certainly exists a maximal weight two Hodge structure on $T$.
Note that the Hodge structure $(T(V), \mathbf C \omega_{V})$ on
$T(V)$ of an abelian surface $V$ is maximal in this sense.
\par
\vskip 4pt
Set $\omega_{i} := \sigma_{i}(\omega)$. Then $(T_{i}, \mathbf C
\omega_{i})$ is a maximal weight two
Hodge structure on $T_{i}$ such that
$\sigma_{i} : (T, \mathbf C \omega) \simeq (T_i, \mathbf C \omega_i)$.
Then, by Theorem 2.1, there is a marked complex $2$-torus
$(A_{i}, \tau_{i})$ such that
$$
\tau_{i} : (H^{2}(A_{i}, \mathbf Z), \mathbf C \omega_{A_{i}})
\simeq (U^{\oplus 3}, \mathbf C \omega_{i})\, .
$$
Then, by the maximality of $T(A_{i})$ and $T_{i}$, one has the
Hodge isometry
$$
\tau_{i} \vert_{T(A_{i})} : (T(A_{i}), \mathbf C \omega_{A_{i}})
\simeq (T_{i}, \mathbf C \omega_{i})
$$
and whence an isometry of the orthogonal lattices
$$
\tau_{i} \vert_{NS(A_{i})} : NS(A_{i}) \simeq S_{i}
$$
for each $i$. Since $S_i$ is hyperbolic, the complex $2$-tori
$A_{i}$ are actually projective, i.e. abelian surfaces.
Since $S_i \not\simeq S_j\; (i\not=j)$ by the choice of $S_{i}$, one
has $NS(A_{i}) \not\simeq NS(A_{j})$ ($i \not = j$) as well.
\par
\vskip 4pt
Set $A := A_{1}$. Then, composing the Hodge isometries
$(\tau \vert_{T(A)}^{-1})$, $\sigma_{i}^{-1}$ and
$\tau_{i} \vert_{T(A_{i})}$, one has a Hodge isometry
$$
(T(A_{i}), \mathbf C \omega_{A_{i}}) \simeq (T(A), \mathbf C \omega_{A})\, .
$$
Therefore $A_{i} \in FM(A)$ by the characterization of $FM(A)$. Thus,
$A_{i} \in \mathcal K(\text{Km}\, A)$ by assertion 1
of Main Theorem. Hence, these $N$ abelian surfaces $A_i$ ($1 \le i \le N$)
satisfy the requirement
in assertion 3 of Main Theorem. \qed

\appendix  % <===== !!!!
%%%%%%%%%%% appendix A %%%%%%%%
\section{}
In this appendix, we shall show the following:
\begin{proposition} \label{prop:AppnA1} Consider the complex $2$-torus
$A = \mathbf C^{2}/L$ defined by a rank $4$ discrete lattice
$L$ in $\mathbf C^2$ which is generated by
$$e_{1} := \Pmatrix  1\\ 0\\ \Endpmatrix\, , \,
e_{2} := \Pmatrix  0\\ 1\\ \Endpmatrix\, , \,
e_{3} := \Pmatrix  \alpha\\ \beta\\ \Endpmatrix\, , \,
e_{4} := \Pmatrix  \gamma\\ \delta\\ \Endpmatrix\, . $$
Assume that $L$ is generic in the sense that $\alpha$, $\beta$, $\gamma$
and $\delta$ are algebraically independent over $\mathbf Q$. Then:
\begin{enumerate}
\item $NS(A) = \{0\}$ and $A$ is not projective.
\item $\hat A \not\simeq A$ and $\hat A \in FM(A)$.
\item $\mathcal K({\rm Km}\, A) = \{ A\}$.
\end{enumerate}
In particular, $\hat A \not\in \mathcal K({\rm Km}\,A)$, i.e.
${\rm Km}\, \hat A \not\simeq {\rm Km}\, A$, while $\hat A \in FM(A)$.
\end{proposition}
This proposition shows certain differences between the Kummer structures
on a non-projective K3 surface and a projective
K3 surface.
\par
\vskip 4pt
\noindent
{\it Proof of  assertion 1.} Let $z_{1}$ and $z_{2}$ be the coordinates of
$\mathbf C^{2}$
corresponding
to the first and the second projections $\mathbf C^{2} \rightarrow \mathbf C$.
We identify $L$ with
$H_{1}(A, \mathbf Z)$.
Let $\langle v_{i} \rangle_{i=1}^{4}$
be the dual basis of $H^{1}(A, \mathbf Z)$ of the basis
$\langle e_{i} \rangle_{i=1}^{4}$ of $H_{1}(A, \mathbf Z)$.
Then, by the description of $e_{i}$, we have
$$dz_{1} = v_{1} + \alpha v_{3} + \gamma v_{4}\, , \, dz_{2} = v_{2} +
\beta v_{3} + \delta v_{4}\, $$
in $H^{1}(A, \mathbf C)$. From this, we obtain
$$\omega := dz_{1} \wedge dz_{2} = v_{1} \wedge v_{2} + \beta v_{1}
\wedge v_{3} + \delta v_{1} \wedge v_{4}
- \alpha v_{2} \wedge v_{3} - \gamma v_{2} \wedge v_{4} +
(\alpha \delta - \beta \gamma)v_{3} \wedge v_{4}\, .$$
Note that $NS(A) = \{\eta \in H^{2}(A, \mathbf Z) \vert (\eta, \omega) = 0\}$.
This due to
the Lefschetz $(1,1)$ Theorem.
Let $\eta \in NS(A)$ and write $\eta$ as
$$\eta = a v_{1} \wedge v_{2} + b v_{1} \wedge v_{3} + c v_{1} \wedge v_{4} +
d v_{2} \wedge v_{3} + e v_{2} \wedge v_{4} +  f v_{3} \wedge v_{4}\, .$$
Here the coefficients $a$, $b$, $\cdots$, $f$ are integers.
By $(\eta, \omega) = 0$, one has
$$a(\alpha \delta - \beta \gamma) + b\gamma - c\alpha + d\delta - e\beta + f
= 0\, .$$
Since $\alpha$, $\beta$, $\gamma$, $\delta$ are algebraically independent
over $\mathbf Q$, we have then
$$a = b = c = d = e = f = 0\, , \,{\rm i.e.}\, , \, \eta = 0\, .$$
Hence $NS(A) = \{0\}$. This also implies that $A$ is not projective. \qed
\par
\vskip 4pt
\noindent
{\it Proof of  assertion 2.} The second statement $\hat A \in FM(A)$
follows from [Mu1],
that the Poincar\'e bundle
induces the equivalence $D(\hat A) \simeq D(A)$.
\par
\vskip 4pt
Let us show the first statement. Let $\hat L$ be the rank 4 discrete
lattice in $\mathbf C^{2}$ generated by
$$f_{1} := \Pmatrix  1\\ 0\\ \Endpmatrix\, , \,
f_{2} := \Pmatrix  0\\ 1\\ \Endpmatrix\, , \,
f_{3} := \Pmatrix  \alpha\\ \gamma\\ \Endpmatrix\, , \,
f_{4} := \Pmatrix  \beta\\ \delta\\ \Endpmatrix\, . $$
Then, by the formula [Sh2, (2.13)], we have $\hat A = \mathbf C^{2}/\hat L$.
Assume that
$A \simeq \hat A$. Then, there is $B \in \text{GL}(2, \mathbf C)$ such that
$Be_{1}\, ,\, Be_{2}\, ,\, Be_{3}\, ,\, Be_{4}\, \in \hat L$.
\par
\vskip 4pt
By $Be_{1}\, ,\, Be_{2} \in \hat L$, there are integers
$m_{i}$, $n_{i}$ such that
$$Be_{1} =
\Pmatrix m_{1} + m_{3}\alpha + m_{4}\beta
\\ m_{2} + m_{3}\gamma + m_{4}\delta \\ \Endpmatrix \, ,$$
$$Be_{2} =
\Pmatrix n_{1} + n_{3}\alpha + n_{4}\beta \\ n_{2} + n_{3}\gamma +
n_{4}\delta \\ \Endpmatrix \, .$$
Then, we have
$$B = \Pmatrix  m_{1} + m_{3}\alpha + m_{4}\beta &
n_{1} + n_{3}\alpha + n_{4}\beta \\
m_{2} + m_{3}\gamma + m_{4}\delta & n_{2} + n_{3}\gamma + n_{4}\delta \\
\Endpmatrix\, .$$
Therefore
$$Be_{3} = \Pmatrix m_{1}\alpha  + m_{3}\alpha^{2} + m_{4}\alpha\beta +
n_{1}\beta  + n_{3}\alpha\beta + n_{4}\beta^{2} \\ m_{2}\alpha  + m_{3}
\alpha\gamma + m_{4}\alpha\delta +
n_{2}\beta  + n_{3}\beta\gamma + n_{4}\beta\delta \\ \Endpmatrix \, .$$
Since $Be_{3} \in \hat L$, there are integers $l_{i}$ such that
$$\Pmatrix m_{1}\alpha  + m_{3}\alpha^{2} + m_{4}\alpha\beta +
n_{1}\beta  + n_{3}\alpha\beta + n_{4}\beta^{2} \\ m_{2}\alpha  + m_{3}
\alpha\gamma + m_{4}\alpha\delta +
n_{2}\beta  + n_{3}\beta\gamma + n_{4}\beta\delta \\ \Endpmatrix
= \Pmatrix l_{1} + l_{3}\alpha + l_{4}\beta \\ l_{2} + l_{3}\gamma + l_{4}
\delta \\ \Endpmatrix \, .$$
Since $\alpha$, $\beta$, $\gamma$, and $\delta$ are algebraically
independent over $\mathbf Q$, this implies that
$$m_{2} = m_{3} = m_{4} = n_{2} = n_{3} = n_{4} = 0\, .$$
However, then $B$ is of the form,
$$B = \Pmatrix  a & b \\ 0 & 0\\ \Endpmatrix\, ,$$
and therefore $\text{det}\, B = 0$, a contradiction.
Hence $A \not\simeq \hat A$. \qed
\par
\vskip 4pt
\noindent
{\it Proof of  assertion 3.} Let $B$ be a complex $2$-torus
such that $\text{Km}\, B \simeq \text{Km}\, A$.
Let us choose an isomorphism $\varphi : \text{Km}\, A \simeq \text{Km}\, B$.
Since $NS(A) = \{0\}$, there is no effective curves on $A$. Therefore
$\text{Km}\, A$ contains no effective curves
other than the $16$ exceptional curves $C_{i}$ arising from the resolution
$\text{Km}\, A \rightarrow A/\langle -1 \rangle$.
Therefore $\text{Km}\, B$ also contains only 16 $\mathbf P^{1}$s, say $D_{i}$,
by $\text{Km}\, B \simeq \text{Km}\, A$.  These $D_{i}$ must then
coincide with
the $16$ exceptional curves arising from the resolution
$\text{Km}\, B \rightarrow B/\langle -1 \rangle$ and satisfy
$\varphi(\cup_{i=1}^{16} C_{i}) = \cup_{i=1}^{16} D_{i}$. Now recall that $A$ can
be reconstructed from $\text{Km}\, A$
as follows: first take the (unique) double cover of
$\text{Km}\, A$ branched along the divisor $\sum_{i=1}^{16}C_{i}$ and
then contract the proper transform of
$\sum_{i=1}^{16}C_{i}$. Likewise $B$ from $\text{Km}\, B$ and
$\sum_{i=1}^{16}D_{i}$ .
Since $\varphi(\sum_{i=1}^{16} C_{i}) = \sum_{i=1}^{16} D_{i}$, it follows
that $\varphi$ then lifts to an isomorphism $\tilde{\varphi} : A \simeq B$.
This completes the proof. \qed
\par
\vskip 4pt
\noindent
{\bf Remark.} In [Ni1, main result], there is a bijective
correspondence
between the set of Kummer structures $\mathcal K(X)$ on $X$ and
$\mathcal C(X)/\text{Aut}(X)$,
where $\mathcal C(X)$ is the set of the configurations of $16$ disjoint
$\mathbf P^{1}$s on $X$. However
$\mathcal C(X)$ is an infinite set and the natural action of
$\text{Aut}(X)$ on $\mathcal C(X)$
is not transitive in general, when $X$ is a projective K3 surface.
It seems much harder to prove our
main theorem by such a geometric method using this
bijective correspondence.


\begin{thebibliography}{[HLOY]}

\bibitem[BM]{BM}
T. Bridgeland, A. Maciocia,
\textit{Complex surfaces with equivalent derived categories},
Math. Zeit. {\bf 236} (2001) 677--697.

\bibitem[BO]{BO}
A. Bondal, D. Orlov,
\textit{Reconstruction of a variety from the derived category and
groups of autoequivalences}, math.AG/9712029, Compositio Math. {\bf 125}
(2001) 327--344.

\bibitem[Cs]{Cs}
J. W. S. Cassels,
\textit{Rational quadratic forms}, Academic Press (1978).

\bibitem[CL]{CL}
H. Cohen, H.W. Lenstra, Jr.,
\textit{Heuristics on class groups of number field},
in Number Theory, Noordwijikerhout 1983,
LNM {\bf 1068} (1984) 33--62.

\bibitem[GH]{GH}
V. Gritsenko, K. Hulek,
\textit{Minimal Siegel modular threefolds}, math.AG/9506017,
Math. Proc. Camb. Phil. Soc. {\bf 123} (1998) 461--485.

\bibitem[HS]{HS}
K. Hulek, G.K. Sankaran,
\textit{The geometry of Siegel modular varieties}, math.AG/9810153,
to appear in Adv. Studies in Pure Math..

%\bibitem[HLOY1]{HLOY1}
%S. Hosono, B.H. Lian, K. Oguiso, S.T. Yau,
%\textit{Autoequivalences of derived category of a K3 surface and monodromy
%transformations}, math.AG/0201047.

\bibitem[HLOY]{HLOY}
S. Hosono, B.H. Lian, K. Oguiso, S.T. Yau,
\textit{Counting Fourier-Mukai partners and applications},
math.AG/0202014.

\bibitem[Iw]{Iw}
H. Iwaniec,
\textit{Almost-primes represented by quadratic polynomials},
Invent. Math. {\bf 47} (1978) 171--188.

\bibitem[Ka]{Ka}
Y. Kawamata,
\textit{D-equivalence and K-equivalence},
math.AG/0205287.

\bibitem[Mo]{Mo}
D. Morrison,
\textit{On K3 surfaces with large Picard number},
Invent. Math. {\bf 75} (1984) 105--121.

\bibitem[Mu1]{Mu1}
S. Mukai,
\textit{Duality between $\mathbf D(X)$ and $\mathbf D(\hat X)$ with
its application to Picard sheaves},
Nagoya. Math. J. {\bf 81} (1981) 101--116.

\bibitem[Mu2]{Mu2}
S. Mukai,
\textit{On the moduli space of bundles on K3 surfaces I},
in ``Vector bundles on algebraic varieties'', Oxford Univ. Press (1987)
341 -- 413.

\bibitem[Na]{Na}
Yo. Namikawa,
\textit{Counter-example to global Torelli problem for irreducible
symplectic manifolds}, math.AG/0110114.

\bibitem[Ni1]{Ni1}
V. V. Nikulin,
\textit{On Kummer surfaces},
Math. USSR. Izv. {\bf 9} (1975) 261--275.

\bibitem[Ni2]{Ni2}
V. V. Nikulin,
\textit{Integral symmetric bilinear forms and some of their geometric
applications},
Math. USSR. Izv. {\bf 14} (1980) 103--167.

\bibitem[Og]{Og}
K. Oguiso,
\textit{K3 surfaces via almost-primes},
math.AG/0110282, Math. Res. Lett. {\bf 9} (2002) 47--63.

\bibitem[Or]{Or}
D. Orlov,
\textit{Equivalences of derived categories and K3 surfaces},
math.AG/9606006, J. Math. Sci. {\bf 84} (1997) 1361--1381.

\bibitem[Sh1]{Sh1}
T. Shioda,
\textit{Some remarks on abelian varieties},
J. Fac. Sci. Univ. Tokyo, Sect IA {\bf 24} (1977) 11--21.

\bibitem[Sh2]{Sh2}
T. Shioda,
\textit{The period map of abelian surfaces},
J. Fac. Sci. Univ. Tokyo, Sect IA {\bf 25} (1978) 47--59.

\bibitem[SM]{SM}
T. Shioda, N. Mitani,
\textit{Singular abelian surfaces and binary quadratic forms},
in Classification of algebraic varieties and complex manifolds,
LNM {\bf 412} (1974) 259--287.

\bibitem[Za]{Za}
D. Zagier
\textit{Zetafunktionen und quadratische Korper : eine Einfuhrung in die
hohere Zahlentheorie},
Springer-Verlag (1981).

\end{thebibliography}
\end{document}